\documentclass[11pt]{article}
\usepackage{amsmath, amsthm, amssymb}
\begin{document}
\title{Caculus of Variation and the $L^{2}$-Bergman Metric on Teichm\"{u}ller Space}

\author{Zheng Huang}
\date{}
\newtheorem{theorem}{Theorem}[section]

\newtheorem{pro}[theorem]{Proposition}
\newtheorem{cor}[theorem]{Corollary}
\newtheorem{lem}[theorem]{Lemma}
\newtheorem{rem}[theorem]{Remark}

\newcommand{\WP}{Weil-Petersson}
\newcommand{\TS}{Teichm\"{u}ller space}
\newcommand{\ms}{moduli space}
\newcommand{\cs}{conformal structure}
\newcommand{\Bm}{Bergman metric}
\newcommand{\cm}{canonical metric}
\newcommand{\hq}{holomorphic quadratic}
\newcommand{\RS}{Riemann surface}
\newcommand{\hm}{harmonic map}
\newcommand{\hym}{hyperbolic metric}
\newcommand{\cd}{complex dimension}
\newcommand{\Bd}{Beltrami differential}
\newcommand{\ts}{tangent space}
\newcommand{\im}{identity map}
\newcommand{\Gc}{Gaussian curvature}
\newcommand{\Ad}{Abelian differential}
\newcommand{\qd}{quadratic differential}
\newcommand{\Hd}{Hopf differential}

\maketitle
\begin{abstract}
The {\cm} on a surface is of nonpositive curvature, so it is natural to study 
{\hm}s between {\cm}s on a surface in a fixed homotopy class. Through this 
approach, we establish the $L^2$-{\Bm} on {\TS} as the second variation of 
energy functionals of chosen families of {\hm}s.   
\end{abstract}
\tableofcontents
\footnotetext{Zheng Huang: Department of Mathematics, 
University of Michigan, Ann Arbor, MI 48109, USA. Email: zhengh@umich.edu}
\footnotetext{AMS subject classification 30F60, 32G15}
\section {Introduction}
We present a geometric analytic approach to the $L^{2}$-{\Bm} on {\TS} of {\RS}s 
in this paper.

{\TS} ${\mathcal {T}}_{g}$ is the space of {\cs}s on  a compact, 
smooth, oriented, closed {\RS} $\Sigma$ of genus $g \ge 1$, where two {\cs}s 
$\sigma$ and $\rho$ are equivalent if there is a biholomorphic map between 
$(\Sigma,\sigma)$ and $(\Sigma,\rho)$ in the homotopy class of the {\im}. When 
$g \ge 2$,  {\TS} ${\mathcal {T}}_{g}$ is naturally a complex manifold of {\cd} 
$3g - 3 > 1$,  and the co{\ts} at $\Sigma$ is identified with $QD(\Sigma)$, 
the space of {\hq} differentials. (\cite {Ah}).

Since ${\mathcal {T}}_{g}$ is a complex manifold, it is natural to study its metric geometry. 
There are several interesting metrics defined on {\TS}, all have advantages and disadvantages. 
These metrics reflect different perspectives of {\TS}. Among those metrics, two are named 
after S. Bergman. One of which, we still call it the {\Bm}, comes from the Bergman kernel 
function of a complex manifold: as a bounded complex domain ${\mathcal {T}}_{g}$ 
carries an invariant K\"{a}hlerian {\Bm}, defined by the line element
\begin{center}
$ds^2 = \sum {\frac{\partial^2 logK(z,z)}{\partial z \partial\bar{z}}} dz_i 
\wedge d{\bar{z}}_j$,
\end{center}
where $K(z,\zeta)$ is the Bergman kernel. This {\Bm} is complete (\cite {Hn}).

The main object of this paper is the other metric which sometimes bears the name of 
Bergman. It is a {\WP} type metric on {\TS}, i.e., it is obtained  from duality by a 
$L^2$ inner product. In order to distinguish from the first {\Bm}, it may be appropriate to 
call this metric the $L^2$-{\Bm}.

From the classical {\RS} theory, the period map $p: \Sigma \rightarrow J_{\Sigma}$ 
embeds the surface $\Sigma$ to its Jacobian $J_{\Sigma}$. The pullback metric of the flat 
metric on $J_{\Sigma}$ via this period map thus defines the so-called {\it {\cm}} or 
{\it {\Bm}} on $\Sigma$, denoted by $\rho_B$. This metric $\rho_B$ is of nonpositive 
{\Gc}, and when $g \ge 2$, the curvature vanishes if and only if the surface is hyperelliptic 
and only at $2g+2$ Weierstrass points (\cite {GR} \cite {L}), in other words, the {\Gc}s 
characterize hyperelliptic surfaces. There is a unique {\cm} in every conformal structure.

The induced $L^2$-Bergman cometric is defined on $QD(\Sigma)$ by $L^2$-norm
\begin{center}
$\| \phi \|_{B}^2 = \int_{\Sigma} \frac{|\phi|^2}{\rho_B}$, 
\end{center}
thus we obtain a metric on {\TS} by duality. This is a Riemannian, 
Hermitian metric, invariant under the mapping class group.

This metric has been studied by Haberman and Jost who showed that it is 
incomplete (\cite {HJ}). Roughly speaking, with respect to the $L^2$-{\Bm}, 
boundary points of the {\ms} $\mathcal{M}_g$ corresponding to pinching a 
nonseparating curve on the surface are at infinite distance from the 
interior, while boundary points of $\mathcal{M}_g$ corresponding to 
pinching a separating curve on the surface are at finite distance from the 
interior. In a sense, the $L^2$-{\Bm} detects topology of the surface.

One of the motivations of this study is to compare the $L^2$-{\Bm} with more 
intensively studied {\WP} metric. These two metrics are both defined from duality 
from $L^2$ inner products, and they are both incomplete. However, the 
$L^2$-{\Bm} does not depend on the uniformization theorem. The difference 
between {\hym} (constant curvature $-1$) and the {\cm} on the surface results in 
different behavior of the induced $L^2$ metrics on {\TS}. The {\WP} metric is of 
negative curvature (\cite {Tr}, \cite {Wp86}), we are yet to understand the 
curvature properties of the $L^2$-{\Bm}.

In this paper, we take a variational approach to the study of the $L^2$-{\Bm}. To do 
so, we fix a {\cs} $(\sigma,z)$ with conformal coordinates $z$. For each {\cm} $\rho$ 
on the surface $\Sigma$, one obtains a {\qd} $\phi(z)dz^{2}$ which is the {\Hd} of the 
unique {\hm} from $\sigma$ to $\rho$. This {\qd} is {\it holomorphic} with respect to 
the {\cs} $(\sigma,z)$, therefore an element of the space $QD(\Sigma)$. We thus obtain 
a map $\phi$ between {\TS} ${\mathcal {T}}_{g}$ and $QD(\Sigma)$, sending $\rho$ to 
$\phi(z)dz^{2}$. We show that this map is a global homeomorphism, hence it provides 
global coordinates to ${\mathcal {T}}_{g}$. The following theorem is an analog to 
Wolf's theorem in the case of {\hym}s (\cite {Wf}).
\begin{theorem}
The map $\phi: {\mathcal {T}}_{g} \rightarrow  QD(\Sigma)$ is a homemorphism.
\end{theorem}

We note that, in the case of {\hym}s, the injectivity of the map $\phi$ is a direct 
application of Bochner's identities and maximum principle, as seen in {\cite {S}}, 
relying on the fact that {\hym} is of constant curvature $-1$. In the case of {\cm}s, 
this is rather difficult since {\cm} has varied curvatures.

With the homeomorphism theorem in mind, we then consider a family of {\hm}s between 
{\cm}s on a surface and show that the second variation of an energy functional is 
the $L^2$-{\Bm} of two infinitesimal cotangent vectors on {\TS}. In the case of 
varying target metrics, we find: 
\begin{theorem}
Let $w(t): (\Sigma, \sigma(z)|dz|^2) \rightarrow 
(\Sigma, \rho(t)|dw|^2)$ be a family of {\hm}s between {\cm}s on surface 
$\Sigma$, where $\rho(0) = \sigma$, for $|t| < \epsilon$ small. Then the 
second variation of the energy functional of $w(t)$, at $t=0$, is given by 
the $L^2$-{\Bm} of infinitesimal {\hq} differentials (up to a constant).
\end{theorem}

Similar result holds in the case of varying domain metrics:
\begin{theorem}
Let $w(s): (\Sigma, \sigma(s)) \rightarrow (\Sigma, \rho)$ be a family of 
{\hm}s between {\cm}s on surface $\Sigma$, where $\sigma(0) = \rho$, for 
$|s| < \epsilon$ small. Then the second variation of the energy
functional of $w(t)$, at $t=0$, is given by the $L^2$-{\Bm} of 
infinitesimal {\hq} differentials (up to a constant).
\end{theorem}

This paper is organized as follows. We introduce the preliminaries in section 2, 
then prove the homeomorphism theorem 1.1 in section 3. Section Four is 
devoted to a variational approach to the study of the $L^2$-{\Bm}, where we 
prove theorem 1.2 (varying the target metric) in $\S 4.1$ and theorem 1.3 
(varying the domain metric) in $\S 4.2$.

The author owes a great debt to, and wishes to thank, Xiaodong Wang and Mike Wolf 
for helpful discussions over the topic. 
\section{Preliminaries}
On a compact {\RS} $\Sigma$ of genus $g>1$, the dimension of the space of 
{\Ad}s of the first kind, or holomorphic one forms, is $g$. There is a natural 
pairing of {\Ad}s defined on this space:
\begin{center}
$<\mu, \nu> = {\frac{\sqrt{-1}}{2}}\int_{\Sigma}\mu\wedge\bar{\nu}$
\end{center}
Let $\{\omega_1, \omega_2, \cdots, \omega_g \}$ be a basis of {\Ad}s, normalized 
with respect to the $A$-cycles of some symplectic homology basis 
$\{A_i, B_i \}_{1 \le i \le g}$, i.e., $\int_{A_i}\omega_j = \delta_{ij}$. 
Thus the period matrix $\Omega_{ij} = \int_{B_i}\omega_j$. One finds that, since 
not all {\Ad}s vanish at the same point according to Riemann-Roch, the period 
matrix is then symmetric with positive definite imaginary part:
$Im \Omega_{ij} = <\omega_i, \omega_j>$ (\cite{FK}).

The {\cm} $\rho_B$ on surface $\Sigma$ is the metric associated to the $(1,1)$ 
form given by 
\begin{center}
$ {\frac{\sqrt{-1}}{2}}\sum_{i,j=1}^{g}
(Im\Omega)_{ij}^{-1}\omega_{i}(z){\bar{\omega}}_{j}(\bar{z})$.
\end{center}

It is not hard to see that this metric is the pull-back of the Euclidean metric 
from the Jacobian variety $J(\Sigma)$ via the period map (\cite {FK}). 
\begin{rem}
It is easy to see that the area of the surface $\Sigma$ with respect to the {\cm} 
is a constant, i.e., $\int_{\Sigma}\rho_B = g$. Sometimes the {\cm} is also 
refered to ${\frac{\rho_B}{g}}$ to unify the surface area. 
\end{rem}

It is known that, when $g \ge 2$, the {\Gc} $K_c$ satisfies $K_c \le 0$ 
(\cite {GR}, \cite {L}), and $K_c(p) = 0$ for some $p \in \Sigma$ if and only if 
$\Sigma$ is hyperelliptic and $p$ is one of the $2g+2$ classical Weierstrass points of 
$\Sigma$ (\cite {L}).

The {\WP} cometric on {\TS} is defined on the space of {\hq} differentials 
$QD(\Sigma)$ by the $L^2$-norm:
\begin{eqnarray}
||\phi||_{WP}^2 = \int_{\Sigma} \frac {|\phi|^2}{\sigma}dzd\bar{z}
\end{eqnarray}
where $\sigma |dz|^2$ is the {\hym} on $\Sigma$. By duality, we obtain a 
Riemannian metric on the tangent space of ${\mathcal{T}}_g$.

The {\it $L^2$-{\Bm}} on ${\mathcal{T}}_g$ is similarly defined by duality from the 
$L^2$-norm
\begin{center}
$\| \phi \|_{B}^2 = \int_{\Sigma} \frac{|\phi|^2}{\rho_B}$. 
\end{center}

We now introduce {\hm}s between {\cm}s on a surface as much of our analysis will 
focus on the techniques of {\hm}s.

For a Lipschitz map $w:(\Sigma, \sigma |dz|^2) \rightarrow 
(\Sigma, \rho |dw|^2)$, where $\sigma |dz|^2$ and $\rho |dw|^2$ are 
metrics on $\Sigma$, and $z$ and $w$ are conformal coordinates on $\Sigma$, 
one follows some notations of Sampson (\cite {S}) to define 
\begin{center}
${\mathcal{H}}(z) = {\frac{\rho(w(z))}{\sigma(z)}}|w_z|^2, {\mathcal{L}}(z) = 
{\frac{\rho(w(z))}{\sigma(z)}}|w_{\bar{z}}|^2$.
\end{center}
Then the energy density of $w$ is simply $e(w)={\mathcal{H}}+{\mathcal{L}}$, and the total energy 
is then given by
\begin{center}
$E(w,\sigma,\rho) = \int_{\Sigma}e\sigma|dz|^2$,
\end{center}
which depends on the target metric and conformal structure of the domain. 
The map $w$ is called {\it harmonic} if it is a critical point of this energy 
functional, i.e., it satisfies Euler-Lagrange equation:
\begin{center}
$w_{z\bar{z}}+ {\frac{\rho_w}{\rho}}w_z w_{z\bar{z}} = 0$.
\end{center}

The $(2,0)$ part of the pullback $w^{*}\rho$ is the so-called {\it {\Hd}}:
\begin{center}
$\phi(z)dz^2 = (w^{*}\rho)^{(2,0)} = \rho w_z {\bar{w}}_zdz^2$.
\end{center}
It is routine to check that $w$ is harmonic if and only if $\phi dz^2 \in QD(\Sigma)$, 
and $w$ is conformal if and only if $\phi=0$.

One also finds that 
\begin{eqnarray}
{\mathcal{H}}(z){\mathcal{L}}(z) = {\frac{\phi {\bar{\phi}}}{\sigma^2}} = {\frac{|\phi|^2}{\sigma^2}} 
\end{eqnarray}
and the Jacobian functional is $J(z) = {\mathcal{H}}(z) - {\mathcal{L}}(z)$.

Now assume both $\sigma$ and $\rho$ are {\cm}s on surface $\Sigma$ (then they 
represent two different {\cs}s unless they are biholomorphic). Since the target 
surface $(\Sigma, \rho)$ has negative Gaussian curvatures almost everywhere, with 
possibly finitely many flat points, the classical theory of {\hm}s guarantees 
that there is a unique {\hm} $w:(\Sigma, \sigma) \rightarrow (\Sigma, \rho)$ in 
the homotopy class of the identity, moreover, this map $w$ is a diffeoemorphism 
with $J > 0$ and ${\mathcal{H}}>0$ (\cite {ES}, \cite {Hr}, \cite {S}, \cite {SY}, 
\cite {CH}, \cite {JS}).

\section{A Homeomorphism}
The method of {\hm}s has been a great computational tool in Teichm\"{u}ller theory (see \cite {J}) . 
In the case of {\hym}s on a compact {\RS}, the second variation of the energy of 
the {\hm} $w=w(\sigma,\rho)$, with respect to the domain metric $\sigma$ (or 
target metric $\rho$) at $\sigma=\rho$, yields the {\WP} metric on 
${\mathcal{T}}_g$ (\cite{Tr2}, \cite{Wf}). In our case of {\cm}s on a {\RS}, we 
prove a homeomorphism theorem, the theorem 1.1, to link {\TS} of {\cm}s to the 
space $QD(\Sigma)$.

To define this map, we fix a point $\sigma$ in {\TS}, with conformal coordinates $z$. 
Thus $\sigma$ is a {\cs} on surface $\Sigma$. For each {\cm} $\rho$ on $\Sigma$, we 
obtain the unique {\hm} $w(\sigma,\rho)$ in the homotopy class of the {\im}, since 
$\rho$ is of nonpositive curvature, with only possibly finitely many flat points on 
$\Sigma$. The associated {\Hd} of the {\hm} $w(\sigma,\rho)$ is then given by 
$\phi(z)dz^2 = \rho(w(z))w_z {\bar{w}}_zdz^2$. Therefore the map 
$\phi: {\mathcal{T}}_g \rightarrow QD(\Sigma)$ which sends $\rho$ to $\phi(z)dz^2$ 
is well defined.

\begin{rem}
Sampson considered this map in the case that of {\hym}s (\cite {S}), and showed that 
it is continuous and one-to-one. Later Wolf showed the map is actually a 
homeomorphism (\cite {Wf}). The condition of constant {\Gc} of {\hym} is essential in 
the argument of proving this map is injective.
\end{rem}

We start with a technical lemma, which is only slightly different than the case of 
{\hym}s, as shown in \cite {Wf}. For {\Hd} $\phi$ corresponding to metric $\rho$, we 
define $\|\phi\| = \int_{\Sigma}|\phi|dzd\bar{z}$. We need to show $\|\phi\|$ is 
approximately the total energy of the {\hm} $w$ in a large scale, i.e.,
\begin{lem}
$\|\phi\| \rightarrow \infty$ if and only if $E(\rho) \rightarrow \infty$
\end{lem}
\begin{proof}
This {\hm} $w(z)$ is naturally quasiconformal, and we write its {\Bd} as 
$\nu = {\frac{w_{\bar{z}}}{w_{z}}}$, and $|\nu|<1$.

We abuse our notation to write $\sigma$ as the domain {\cm}, and 
$dA = \sigma(z)dzd\bar{z}$ is the area element of the domain surface. Recall from 
section two, we have density functions 
${\mathcal{H}}(z) = {\frac{\rho(w(z))}{\sigma(z)}}|w_z|^2$ and 
${\mathcal{L}}(z) = {\frac{\rho(w(z))}{\sigma(z)}}|w_{\bar{z}}|^2$. The total energy 
is $E(\rho)=\int_{\Sigma}({\mathcal{H}}(z)+{\mathcal{L}}(z))dA$, and the Jacobian 
determinant of the map is $J(z) = {\mathcal{H}}(z)- {\mathcal{L}}(z)$. Note that 
${\mathcal{H}}(z) > {\mathcal{L}}(z) \ge 0$.

It is not hard to see that 
\begin{center}
$|\phi|^{2} = \sigma^{2}{\mathcal{H}}{\mathcal{L}}$, and 
$|\nu|^{2} = {\frac{{\mathcal{L}}}{{\mathcal{H}}}} < 1$.
\end{center}
Therefore, we now have 
\begin{eqnarray*}
\|\phi\| &=&  \int_{\Sigma}|\phi|dzd\bar{z} = 
\int_{\Sigma}{\mathcal{H}}|\nu|\sigma dzd\bar{z} \\
&=&\int_{\Sigma}{\mathcal{H}}|\nu| dA < \int_{\Sigma}{\mathcal{H}}dA \\
&\le& \int_{\Sigma}({\mathcal{H}}+{\mathcal{L}})dA = E(\rho).
\end{eqnarray*}
For the opposite direction, we find 
${\mathcal{L}} \le \sqrt{{\mathcal{H}}{\mathcal{L}}}$ since 
${\mathcal{L}}<{\mathcal{H}}$, and therefore, 
\begin{eqnarray*}
E(\rho) &=&   \int_{\Sigma}({\mathcal{H}}+{\mathcal{L}})dA =  \int_{\Sigma}J(z)dA + 
2 \int_{\Sigma}{\mathcal{L}}dA \\
&\le& Area(\Sigma,\sigma) + 2\int_{\Sigma}\sqrt{{\mathcal{H}}{\mathcal{L}}}dA \\
&=& g+2\int_{\Sigma}{\frac{|\phi|}{\sigma}}dA = g+ 2\|\phi\|.
\end{eqnarray*}
Here we used the fact that $Area(\Sigma,\sigma) = g$, as pointed out in remark 2.1. 
This completes the proof of the lemma.
\end{proof}
We now start to prove the homeomorphism theorem.
\begin{proof} (Proof of theorem 1.1):
It is clear that this map is continuous because of the uniqueness of {\hm} $w$ in 
the homotopy class of the {\im}. We want to show this map is a local diffeoemorphism 
and proper.

We firstly notice that the map $\phi$ is a local diffeomorphism. To see this, we 
consider a sufficiently small neighborhood of $\sigma$ and a family of {\hm}s 
$w(t):(\Sigma, \sigma) \rightarrow (\Sigma, \rho(t))$ between {\cm}s near $t=0$, 
where $\rho(t)$ is a family of {\cm}s with $\rho(0) = \sigma$. It is easy to see that 
$w(0) = z$, the {\im}. Associated {\Hd}s of this family are given by 
$\phi(t)dz^2 = \rho(t)w_z(t) {\bar{w}}_z(t)dz^2$ with $\phi(0) = 0$. we take 
$t$-derivative on $\phi(t)$ at $t=0$ to find that
 \begin{center}
 ${\frac{d\phi(t)}{dt}}|_{t=0} = \rho(0)w_{z}(0){\frac{d{\bar{w}_{z}}(t)}{dt}}|_{t=0} 
 = \sigma {\frac{d{\bar{w}_{z}}(t)}{dt}}|_{t=0}$.
 \end{center}
This shows that ${\frac{dw(t)}{dt}}|_{t=0}$ is conformal, provided that 
${\frac{d\phi(t)}{dt}}|_{t=0} = 0$. So the map $d\phi$ is nonsingular, and $\phi$ is a 
local diffeomorphism by applying inverse function theorem.

We then apply a slightly rearranged argument of Wolf (\cite{Wf}) (on {\hym}s) to show 
map $\phi$ is proper.

Given that ${\mathcal{T}}_{g}$ and $QD(\Sigma)$ are finite dimensional spaces, and from 
lemma 3.2, it suffices to show the energy function $E(\rho)$ is a proper map from {\TS} 
to $\Re$ (see theorem 2.7.1, \cite{B}). In other words, we need to show the set 
$B=\{\rho \in {\mathcal{T}}_{g}: E(\rho) \le K \}$ is a compact subset of 
${\mathcal{T}}_{g}$. Without loss of generality, we assume 
$id:(\Sigma, \sigma) \rightarrow (\Sigma,\rho)$ is harmonic, or we can choose $w^{*}\rho$ 
to represent the equivalency class $[\rho]$.

Consider a geodesic ball $B(x_{0},\delta)$ for some $x_{0}$ in domain surface 
$(\Sigma, \sigma)$, where positive constant $\delta < min\{1,inj_{\sigma}(\Sigma)^{2}\}$, 
where $inj_{\sigma}(\Sigma)$ is the injectivity radius of $\Sigma$ with respect to the metric 
$\sigma$. Notice that a {\hm} between surfaces does not depend on the choice of metrics on 
the domain, but on the choice of {\cs}s of the domain surface. Therefore, we can choose 
$\sigma$ to be hyperbolic in this argument, and then introduce polar coordinates 
$(r,\theta)$ in the hyperbolic disk $B(x_{0},\delta)$ 
so that $\sigma = dr^{2}+sinh^{2}(r)d\theta^{2}$.

For $r < \sqrt{\delta}<1$, we have $sinh(r)<2r$ and then
\begin{center}
$\int_{\delta}^{\sqrt{\delta}}{\frac{dr}{sinh(r)}} > {\frac{1}{2}}\int_{\delta}^{\sqrt{\delta}}{\frac{dr}{r}}  
= {\frac{1}{4}}|log\delta|$.
\end{center}

Now considering the annulus $A(x_{0})=A(x_{0},\delta,\sqrt{\delta})$ centered at $x_{0}$ of inner and outer radii $\delta$ and 
$\sqrt{\delta}$, respectively, in domain metric $\sigma$, we apply the upper bound of the energy to find
\begin{eqnarray*}
\int \!\!\! \int_{A(x_{0})}\|{\frac{\partial}{\partial \theta}}\|_{\rho}^{2}{\frac{dr}{sinh(r)}}d\theta 
& \le & \int \!\!\! \int_{A(x_{0})}\|{\frac{\partial}{\partial r}}\|_{\rho}^{2} + 
{\frac{1}{sinh^{2}(r)}}\|{\frac{\partial}{\partial \theta}}\|_{\rho}^{2}sinh(r)drd\theta \\
& \le & 2\int_{\Sigma}e(z)dA = 2E(\rho) \le 2K.
\end{eqnarray*}
Thus there exists $\delta < r < \sqrt{\delta}$ such that 
\begin{center}
$\int_{0}^{2\pi}\|{\frac{\partial}{\partial \theta}}\|_{\rho}^{2}d\theta \le {\frac{8K}{|log\delta |}}$.
\end{center}
For this $r$ and two points $x_{3}$ and $x_{4}$ on the boundary of the disk $B_{\sigma}(x_{0},r)$, and two points $x_{1}$ 
and $x_{2}$ in $B_{\sigma}(x_{0},\delta)$, we now have, 
\begin{eqnarray*}
d_{\rho}(w(x_{1}),w(x_{2})) &=&  d_{\rho}(x_{1},x_{2}) \le d_{\rho}(x_{3},x_{4}) \\
&\le& \int_{0}^{2\pi}\|{\frac{\partial}{\partial \theta}}\|_{\rho}d\theta \le 
2\pi\sqrt{\int_{0}^{2\pi}\|{\frac{\partial}{\partial \theta}}\|_{\rho}^{2}d\theta} \\
&\le&{\frac{4\sqrt{2K}\pi}{|log\delta |}}.
\end{eqnarray*}
Applying the Courant-Lebesgue lemma, we conclude that the energy $E(\rho)$ is proper, 
hence so is map $\phi$.

We have showed that the map $\phi$ is a proper local diffeomorphism between {\TS} and 
$QD(\Sigma)$. It is clear that ${\mathcal{T}}_{g}$ is path-connected and $QD(\Sigma)$ is a 
simply connected Hausdorff space. Standard theory of covering map between manifolds implies 
such a local homeomorphism is actually a global homeomorphism. This completes the proof 
of theorem 1.1.
\end{proof}
\section{Variations of the Energy}
In the previous section, we showed that the map $\phi: {\mathcal{T}}_g \rightarrow QD(\Sigma)$ 
is a homeomorphism, thus its inverse map $\phi^{-1}$ provides coordinates for any {\cm} 
$\rho \in {\mathcal{T}}_g$. We will study these coordinates in this section, i.e., we apply 
variational approach to derive the infinitesimal $L^2$-Bergman norm on {\TS}. We will separate 
the cases where either target {\cm}s are varying or domain {\cm}s are varying, in subsections 4.1 
and 4.2, respectively.

\subsection{Varying the Target}
To prove theorem 1.2, we need to develop some infinitesimal calculations for the variation of a {\hm} 
between {\cm}s. This technique is an analog to that of Wolf's on the case of {\hym}s, which plays an 
important role in studying the {\WP} geometry of {\TS}.

Now we consider a family of {\hm}s $w(t):(\Sigma, \sigma) \rightarrow (\Sigma, \rho(t))$ 
between {\cm}s, where $w(t)$ varies real analyticaly in $t$ for $|t| < \epsilon$, 
and $\rho(t)$ is a family of {\cm}s with $\rho(0) = \sigma$, therefore $w(0) = z$. 
Associated {\Hd}s are given by $\phi(t)dz^2 = \rho(t)w_z(t) {\bar{w}}_z(t)dz^2$ 
with $\phi(0) = 0$.

For $t=(t^{\alpha}, t^{\beta})$, denote 
$\phi_{\alpha} = {\frac{d\phi(t)}{dt^{\alpha}}}|_{t=0}$ and 
$\phi_{\beta} = {\frac{d\phi(t)}{dt^{\beta}}}|_{t=0}$ as infinitesimal {\hq} differentials.

Recall that the holomorphic and antiholomorphic functions of this family of {\hm}s are 
\begin{center}
${\mathcal{H}}(t) = {\frac{\rho(w(t))}{\sigma(z)}}|w_z(t)|^2, {\mathcal{L}}(t) = 
{\frac{\rho(w(t))}{\sigma(z)}}|w_{\bar{z}}(t)|^2$.
\end{center}
We denote ${\mathcal{H}}_{\alpha} = {\frac{d{\mathcal{H}}(t)}{dt^{\alpha}}}|_{t=0}$ 
and ${\mathcal{L}}_{\alpha} = {\frac{d{\mathcal{L}}(t)}{dt^{\alpha}}}|_{t=0}$, also 
${\mathcal{L}}_{\alpha \bar{\beta}} = {\frac{d^2}{dt^{\alpha}\bar{dt^{\beta}}}}|_{t=0}{\mathcal{L}}(t)$, and 
we assign similiar meaning for ${\mathcal{H}}_{\alpha \bar{\beta}}$ and $E_{\alpha \bar{\beta}}$.

We also write $K(t) =K(\rho(t)) = -{\frac{1}{2}}\Delta_{\rho}log\rho$ as the {\Gc} of the 
metric $\rho(t)$, and assign obvious meaning to $K_{\alpha}$ and $K_{\alpha \bar{\beta}}$. 
Since $K(\sigma) \le 0$ and is negative everywhere except possibly finitely many points, 
it is not hard to see that the operator $\Delta_{\sigma}+2K(\sigma)$ is invertible on 
$(\Sigma,\sigma)$, and we denote $D_B = -2(\Delta_{\sigma}+2K(\sigma))^{-1}$.
\begin{lem}
For this family of {\hm}s $w(t)$, the following holds:
\begin{itemize}
\item[(i)] 
${\mathcal{H}}(0) = 1$ and ${\mathcal{L}}(0)=0$;
\item[(ii)] 
${\mathcal{L}}_{\alpha} \equiv 0$, ${\mathcal{H}}_{\alpha} = D_{B}(K_{\alpha})$, 
and $\int_{\Sigma}{\mathcal{H}}_{\alpha}\sigma = 0$;
\item[(iii)]
${\mathcal{L}}_{\alpha \bar{\beta}} = {\frac{\phi_{\alpha}{\bar{\phi}}_{\beta}}{\sigma^2}}$; 
\item[(iv)]
${\mathcal{H}}_{\alpha \bar{\beta}} = D_{B}(K_{\alpha \bar{\beta}}) + 
D_{B}(K_{\bar{\beta}}D_{B}(K_{\alpha})) + D_{B}(K_{\alpha}D_{B}(K_{\bar{\beta}}))$ \\ 
$-D_{B}(K(\sigma){\frac{\phi_{\alpha}{\bar{\phi}}_{\beta}}{\sigma^2}}) 
-{\frac{1}{2}}D_{B}(\Delta_{\sigma}(D_{B}(K_{\alpha})D_{B}(K_{\bar{\beta}})))$.
\end{itemize}
\end{lem}
\begin{proof}
\begin{itemize}
\item[(i)] 
\noindent
This is true since the map $w(t)$ is the {\im} at time $t=0$. 
\item[(ii)] 
Recalling formula (2): 
\begin{center}
${\mathcal{H}}(t){\mathcal{L}}(t) = {\frac{\phi(t) {\bar{\phi}}(t)}{\sigma^2}}$, 
\end{center}
\noindent
we take $t$-derivative at $t=0$, to find that
\begin{center}
${\mathcal{H}}_{\alpha}{\mathcal{L}}(0) + {\mathcal{H}}(0){\mathcal{L}}_{\alpha} = 
{\frac{\phi_{\alpha}{\bar{\phi}}(0)+\phi(0){\bar{\phi}}_{\alpha}}{\sigma^2}}$.
\end{center}
\noindent
The righthand side is zero as $\phi(0) = 0$. Therefore ${\mathcal{L}}_{\alpha} \equiv 0$.
\\
\noindent
We notice that 
$\int_{\Sigma}({\mathcal{H}}(t)-{\mathcal{L}}(t))\sigma = \int_{\Sigma}J(t)\sigma = g$ is 
independent of the parameter $t$. Therefore
\begin{center}
$\int_{\Sigma}{\mathcal{H}}_{\alpha}\sigma = \int_{\Sigma}{\mathcal{L}}_{\alpha}\sigma = 0$.
\end{center}
\noindent
From standard Bochner identities, we have 
\begin{eqnarray}
\Delta_{\sigma}log{\mathcal{H}} = 2K(\sigma) - 2K(\rho)({\mathcal{H}}-{\mathcal{L}}).
\end{eqnarray}
\noindent
Therefore 
$\Delta_{\sigma}log{\mathcal{H}}(t)=2K(\sigma)-2K(\rho(t))({\mathcal{H}}(t)-{\mathcal{L}}(t))$ 
\noindent
and 
\begin{eqnarray*}
\Delta_{\sigma}{\mathcal{H}}_{\alpha} &=& \Delta_{\sigma}{\frac{{\mathcal{H}}_{\alpha}}{{\mathcal{H}}(0)}} \\
&=& -2K_{\alpha}({\mathcal{H}}(0)-{\mathcal{L}}(0)) - 2K(\rho(0))({\mathcal{H}}_{\alpha}-{\mathcal{L}}_{\alpha}) \\
&=& -2K_{\alpha} - 2K(\sigma){\mathcal{H}}_{\alpha}.
\end{eqnarray*}
\noindent
We now obtain $(\Delta_{\sigma}+2K(\sigma)){\mathcal{H}}_{\alpha} = -2K_{\alpha}$ and then 
\begin{center}
${\mathcal{H}}_{\alpha} = -2(\Delta_{\sigma}+2K(\sigma))^{-1}(K_{\alpha})=D_{B}(K_{\alpha})$.
\end{center}
\noindent
\item[(iii)]
To calculate next variation, we consider formula 
${\mathcal{H}}(t){\mathcal{L}}(t) = {\frac{\phi(t) {\bar{\phi}}(t)}{\sigma^2}}$ again. 
We find 
\begin{center}
${\mathcal{H}}_{\alpha \bar{\beta}}{\mathcal{L}}(0) + 
{\mathcal{H}}_{\bar\beta}{\mathcal{L}}_{\alpha} + {\mathcal{H}}_{\alpha}{\mathcal{L}}_{\bar\beta} + 
{\mathcal{H}}(0){\mathcal{L}}_{\alpha \bar{\beta}} = {\frac{\phi_{\alpha}{\bar{\phi}}_{\beta}}{\sigma^2}}$.
\end{center}
\noindent
Therefore ${\mathcal{L}}_{\alpha \bar{\beta}} = {\frac{\phi_{\alpha}{\bar{\phi}}_{\beta}}{\sigma^2}}$.
\noindent
\item[(iv)]
We take second $t$-derivative from (3) to find
\begin{eqnarray*}
\Delta_{\sigma}({\mathcal{H}}_{\alpha \bar{\beta}} - {\mathcal{H}}_{\alpha}{\mathcal{H}}_{\bar{\beta}}) &=& 
-2K_{\alpha\bar{\beta}} - 2K_{\alpha}{\mathcal{H}}_{\bar\beta} -2K_{\bar\beta}{\mathcal{H}}_{\alpha}\\
&-& 2K(\sigma)({\mathcal{H}}_{\alpha \bar{\beta}}-{\mathcal{L}}_{\alpha \bar{\beta}}),
\end{eqnarray*}
\noindent
then we obtain that
\begin{eqnarray*}
(\Delta_{\sigma}+2K(\sigma))({\mathcal{H}}_{\alpha \bar{\beta}}) &=& \Delta({\mathcal{H}}_{\alpha}{\mathcal{H}}_{\bar{\beta}}) 
-2K_{\alpha\bar{\beta}} - 2K_{\alpha}{\mathcal{H}}_{\bar\beta} -2K_{\bar\beta}{\mathcal{H}}_{\alpha}\\
&-& 2K(\sigma)({\mathcal{H}}_{\alpha \bar{\beta}}-{\mathcal{L}}_{\alpha \bar{\beta}}).
\end{eqnarray*}
\noindent
Now we apply formulas ${\mathcal{H}}_{\alpha} = D_{B}(K_{\alpha})$, and 
${\mathcal{H}}_{\bar\beta} = D_{B}(K_{\bar\beta})$, and ${\mathcal{L}}_{\alpha \bar{\beta}} = 
{\frac{\phi_{\alpha}{\bar{\phi}}_{\beta}}{\sigma^2}}$ to above equation to complete 
the proof of this lemma.
\end{itemize}
\end{proof}
\begin{rem}
It is very interesting to compare our situation with the case of varations of a {\hm} 
between {hym}s on the surface. If we assume all metrics on the surface are 
hyperbolic with constant {\Gc}, under the same notations, then we have the following 
comparison: (i) and (iii) in lemma 4.1 hold; (ii) also holds except furthermore, 
${\mathcal{H}}_{\alpha} \equiv 0$, i.e., the holomorphic energy reaches its minimum at 
time zero; (iv) of the lemma takes the form of ${\mathcal{H}}_{\alpha \bar{\beta}} = 
D({\frac{\phi_{\alpha}{\bar{\phi}}_{\beta}}{\sigma^2}})$, where 
$D=-2(\Delta_{\sigma}-2)^{-1}$ is a compact, self-adjoint operator. Operator $D_B$ 
in (iv) of lemma 4.1 is not self-adjoint for $L^2$ functions, while $-K(\sigma)D_B$ 
is, and coincides with operator $D$ when $K \equiv -1$.  
\end{rem}
We now consider the variations of the corresponding total energy $E(t)$ of the family 
$w(t)$ near $t=0$, i.e., we show theorem 1.2 in following equivalent form:
\begin{theorem}
${\frac{d^2}{dt^{\alpha}\bar{dt^{\beta}}}}|_{t=0}E(t) = 2<\phi_{\alpha},\phi_{\beta}>_B$.
\end{theorem}
\begin{proof}
The total energy is 
\begin{eqnarray*}
E(t) &=& \int_{\Sigma}({\mathcal{H}}(t) + {\mathcal{L}}(t))\sigma \\
&=& \int_{\Sigma}({\mathcal{H}}(t) - {\mathcal{L}}(t))\sigma + 
2\int_{\Sigma}{\mathcal{L}}(t)\sigma \\
&=& \int_{\Sigma}J(t)\sigma + 2\int_{\Sigma}{\mathcal{L}}(t)\sigma \\
&=& g + 2\int_{\Sigma}{\mathcal{L}}(t)\sigma \ge g,
\end{eqnarray*}
since $E(0) = g$ is equal to the area of the surface. Thus $E(t)$ reaches its 
global minimum $g$ at $t=0$ from (i) of lemma 4.1.

From (ii) of lemma 4.1, it is easy to see that $t=0$ is a critical point of $E(t)$ as 
\begin{center}
$E_{\alpha} = \int_{\Sigma}({\mathcal{H}}_{\alpha}+{\mathcal{L}}_{\alpha})\sigma =0$.
\end{center}
and 
\begin{eqnarray*}
E_{\alpha \bar{\beta}} &=& \int_{\Sigma}({\mathcal{H}}_{\alpha \bar{\beta}}+
{\mathcal{L}}_{\alpha \bar{\beta}})\sigma \\
&=& 2\int_{\Sigma}{\mathcal{L}}_{\alpha \bar{\beta}}\sigma 
= 2\int_{\Sigma}{\frac{\phi_{\alpha}{\bar{\phi}}_{\beta}}{\sigma^2}}\sigma \\
&=& 2<\phi_{\alpha},\phi_{\beta}>_B.
\end{eqnarray*}
\end{proof}
\subsection{Varying the Domain}
In this subsection, we consider a family of {\hm}s between fixed target metric 
and varying domain metrics.

Again, since the target metric is negatively curved (except possibly finitely many 
flat points), we have the existence and uniqueness of a {\hm} in the homotopy 
class of the identity, and this map is a diffeomorphism. In other words, let 
$w(s):(\Sigma, \sigma(s)|dz|^2) \rightarrow 
(\Sigma, \rho |dw|^2)$ be this family 
of {\hm}s near the identity map, where $w(s)$ varies real analytically in $s$ for 
$|s| < \epsilon$, and $\sigma(s)$ is a family of {\cm}s with $\sigma(0) = \rho$, 
therefore $w(0) = z$. Associated {\Hd}s are given by 
$\phi(s)dz(s)^2 = \rho w_z(s){\bar{w}}_z(s)dz(s)^2$ with $\phi(0) = 0$.

For $s=(s^a,s^b)$, similar to last subsection, we denote 
$\phi_a = {\frac{\partial\phi(s)}{\partial s^a}}|_{s=0}$ and 
$\phi_b = {\frac{\partial\phi(s)}{\partial s^b}}|_{s=0}$, and assign similar meanings 
to ${\mathcal{H}}_a$, and ${\mathcal{L}}_{a\bar{b}}$, etc.

 Let $K(s) = -{\frac{1}{2}}\Delta_{\sigma(s)}log\sigma(s)$ be the {\Gc} of the surface 
$(\Sigma, \sigma(s))$ and denote $K_a = {\frac{\partial K(s)}{\partial s^a}}|_{s=0}$. 
Again, since $K(\rho) \le 0$ and is negative everywhere except possibly finitely many 
points, the operator $\Delta_{\rho}+2K(\rho)$ is invertible on $(\Sigma,\rho)$, and we 
denote $D'_B = -2(\Delta_{\rho}+2K(\rho))^{-1}$. This operator $D'_B$ is not 
self-adjoint for $L^2$ functions.

We firstly calculate the variations of these two density functions  ${\mathcal{H}}(s)$ and 
${\mathcal{L}}(s)$. It is interesting to notice the difference with the case of varying the target 
showed in lemma 4.1.
\begin{lem}
For this family of {\hm}s $w(s)$, the following holds:
\begin{itemize}
\item[(i)] 
${\mathcal{H}}(0) = 1$ and ${\mathcal{L}}(0)=0$;
\item[(ii)] 
${\mathcal{L}}_{a} \equiv 0$,  ${\mathcal{H}}_a = -D'_B(K_a)$, 
and $\int_{\Sigma}{\mathcal{H}}_{a}\sigma = 0$;
\item[(iii)]
${\mathcal{L}}_{a \bar{b}} = {\frac{\phi_{a}{\bar{\phi}}_{b}}{\sigma^2}}$.
\end{itemize}
\end{lem}
\begin{proof}
\begin{itemize}
\item[(i)] 
\noindent
It is true since the map $w(0)$ is the {\im}. 
\item[(ii)] 
We take $s^a$-derivative of ${\mathcal{H}}(s){\mathcal{L}}(s)={\frac{\phi(s)\bar{\phi}(s)}{\sigma^2(s)}}$ to find
\begin{center}
${\mathcal{H}}_a {\mathcal{L}}(0) + {\mathcal{H}}(0){\mathcal{L}}_a  = {\frac{\phi_a \bar\phi(0) + \phi(0)\bar\phi_a}{\sigma^2 (0)}} +  
[{\frac{\partial ({\frac{1}{\sigma^2 (s)}})}{\partial s^a}}]|_{s=0}|\phi(0)|^2$,
\end{center}
and this implies ${\mathcal{L}}_a=0$, for $\phi(0)=0$.

Therefore $\int_{\Sigma}{\mathcal{H}}_{a}\sigma = \int_{\Sigma}{\mathcal{L}}_{a}\sigma = 0$.

To calculate ${\mathcal{H}}_{a}$, recalling formula (3):
\begin{center}
$\Delta_{\sigma(s)}log{\mathcal{H}}(s) = 2K(\sigma(s)) - 2K(\rho)({\mathcal{H}}(s)-{\mathcal{L}}(s))$,
\end{center}
we find that 
\begin{eqnarray*}
{\frac{\partial}{\partial s^a}}|_{s=0}(\Delta_{\sigma(s)})log{\mathcal{H}}(0) + \Delta_{\rho}({\mathcal{H}}_a) 
&=& 2K_a - 2K(\rho)({\mathcal{H}}_a-{\mathcal{L}}_a) \\
&=& 2K_a - 2K(\rho)({\mathcal{H}}_a).
\end{eqnarray*} 
Therefore $(\Delta_{\rho}+2K(\rho))({\mathcal{H}}_a) = 2K_a$, and so ${\mathcal{H}}_a = -D'_B(K_a)$.
\item[(iii)] 
For the second variation of ${\mathcal{L}}(s)$, we have 
\begin{eqnarray*}
{\mathcal{L}}_{a \bar{b}} &=& {\mathcal{H}}_{a \bar{b}}{\mathcal{L}}(0) + {\mathcal{H}}_{\bar{b}}{\mathcal{L}}_a + {\mathcal{H}}_a {\mathcal{L}}_{\bar b} 
+ {\mathcal{H}}(0){\mathcal{L}}_{a \bar{b}} \\
&=& {\frac{\phi_a \bar\phi_b}{\rho^2}} + 
[{\frac{\partial ({\frac{1}{\sigma^2 (s)}})}{\partial {\bar{s}^b}}}]|_{s=0}[\phi_a \bar\phi(0) + \phi(0)\bar\phi_a]\\
&+& [{\frac{\partial ({\frac{1}{\sigma^2 (s)}})}{\partial s^a}}]|_{s=0}[\phi_b \bar\phi(0) + \phi(0)\bar\phi_b] 
+ [{\frac{\partial^2 ({\frac{1}{\sigma^2 (s)}})}{\partial s^a \partial \bar{s}^b}}]|_{s=0}[\phi(0)\bar\phi(0)] \\
&=&{\frac{\phi_a \bar\phi_b}{\rho^2}}.
\end{eqnarray*} 
\end{itemize}
\end{proof}

Now we show the equivalent form of theorem 1.3:
\begin{theorem}
${\frac{\partial^2}{\partial s^a \bar{\partial s^b}}}|_{s=0}E(s) = 2<\phi_a,\phi_b>_B$.
\end{theorem}
\begin{proof}
The total energy is now
\begin{eqnarray*}
E(s) &=& \int_{\Sigma}({\mathcal{H}}(s) + {\mathcal{L}}(s))\sigma(s) \\
&=& \int_{\Sigma}J(s)\sigma(s) + 2\int_{\Sigma}{\mathcal{L}}(s)\sigma(s) \\
&=& g + 2\int_{\Sigma}{\mathcal{L}}(s)\sigma(s) \ge g,
\end{eqnarray*}
where $E(0) = g$ reaches the global minimum.

Together with ${\mathcal{L}}(0) =0$, we then find 
\begin{center}
$E_a = 2 {\frac{\partial}{\partial s^a}}|_{s=0}[\int_{\Sigma}{\mathcal{L}}(s)\sigma(s)] = 0$,
\end{center}
and then $s=0$ is also a critical point of $E(s)$.

Now we consider ${\frac{\partial^2}{\partial s^a \bar{\partial s^b}}}|_{s=0}E(s)$ from lemma 4.4. 
We apply ${\mathcal{L}}(0) = {\mathcal{L}}_a = {\mathcal{L}}_{\bar{b}} = 0$ to find 
\begin{eqnarray*}
{\frac{\partial^2}{\partial s^a \bar{\partial s^b}}}|_{s=0}E(s) &=& 
2{\frac{\partial^2}{\partial s^a \partial \bar{s}^b}}|_{s=0}\{\int_{\Sigma}{\mathcal{L}}(s)\sigma(s)\} \\
&=& 2\int_{\Sigma}{\frac{\phi_a \bar\phi_b}{\rho^2}} \rho \\
&=& 2<\phi_a,\phi_b>_B.
\end{eqnarray*} 
\end{proof}
\begin{rem}
{\TS} is a complex manifold (when $g \ge 2$), so it has its own complex structure. For Riemannian metrics on this 
complex manifold, it is ideal that metrics are compatible with the complex structure. Ahlfors (\cite {Ah}) showed that the 
{\WP} metric is K\"{a}hlerian. From the definition, we know that the $L^{2}$-{\Bm} is an Hermitian metric, yet it is 
unknown if the $L^{2}$-{\Bm} is actually K\"{a}hlerian. 
\end{rem}

\end{document}